\UseRawInputEncoding
\documentclass[11pt,makeidx]{amsart}

\usepackage{amscd}
\usepackage{amsmath}
\usepackage{amsfonts}
\usepackage{amssymb}

\usepackage{color}
\definecolor{dark_purple}{rgb}{0.4, 0.0, 0.4}

\usepackage{hyperref}
\hypersetup{
    colorlinks=true,       
    linkcolor= red,          
    citecolor=blue,        
    filecolor=magenta,      
    urlcolor=cyan           
}

\input xy
\xyoption{all}

\def \qed {\hbox{}\nobreak
\vrule width 1.4mm height 1.4mm depth 0mm \par \goodbreak \smallskip}

\font\twelveroman=cmr10 at 12pt

\font\sevenroman=cmr7
\font\fiveroman=cmr5
\newfam\romanfam
\textfont\romanfam=\twelveroman
\scriptfont\romanfam=\sevenroman
\scriptscriptfont\romanfam=\fiveroman

\newcommand{\Z}{{{\mathbf Z}}}

\newcommand{\C}{{{\mathbf C}}}
\newcommand{\Q}{{{\mathbf Q}}}
\newcommand{\A}{{{\mathbf A}}}

\newcommand{\Zp}{{\Z_p}}

\newcommand{\Qp}{{\Q_p}}

\newcommand{\h}{{\mathbf h}}

\newcommand{\T}{{\mathbf T}}

\DeclareMathOperator{\Frob}{{\mathrm Frob}}

\newcommand{\CL}{{\mathcal L}}

\newcommand{\CH}{{{\mathcal H}}}

\newcommand{\CF}{{\mathcal F}}

\newcommand{\n}{{\mathfrak n}}

\newcommand{\q}{{\mathfrak q}}
\newcommand{\m}{{\mathfrak m}}

\newcommand{\lb}{\lambda}

\newcommand{{\Art}}{{\mathrm{Art}}}
\newcommand{{\Bcris}}{{\mathrm B_{cris}}}
\newcommand{{\Bst}}{{\mathrm B_{st}}}
\newcommand{{\BdR}}{{\mathrm B_{dR}}}
\newcommand{{\Dcris}}{{\mathrm D_{cris}}}
\newcommand{{\Dst}}{{\mathrm D_{st}}}
\newcommand{{\DdR}}{{\mathrm D_{dR}}}
\newcommand{{\Fil}}{{\mathrm{Fil}}}
\newcommand{{\GQp}}{{G_{\Q_p}}}
\newcommand{{\GQ}}{{G_{\Q}}}
\newcommand{\Gal}{{\mathrm{Gal}}}

\newcommand{\cyc}{\mathrm{cyc}}
\newcommand{\Iw}{\mathrm{Iw}}

\newcommand{\Spec}{\mathrm{Spec}\,}

\newcommand{{\Cusp}}{{{C}}}
\newcommand{{\CuspL}}{{{C}}}

\newcommand{\dis}{\displaystyle}
\newcommand{\ot}{\otimes}

\newcommand{\bQp}{\bQ_p}

\newcommand{\bQ}{{\overline{\Q}}}

\newcommand{\bslh}{\backslash}

\makeatletter
\@addtoreset{equation}{section}

\makeatother

\newtheorem{defi}[equation]{\bf Definition}
\newtheorem{coro}[equation]{\bf Corollary}
\newtheorem{lem}[equation]{\bf Lemma}
\newtheorem{thm}[equation]{\bf Theorem}
\newtheorem{prop}[equation]{\bf Proposition}
\newtheorem{conj}[equation]{\bf Conjecture}

\theoremstyle{definition}

\author{Eric Urban}
\address{
Department of Mathematics\\
 Columbia University\\
 2990 Broadway\\
 New York, NY 10027
 }
 \email{ urban@math.columbia.edu}

\title{On Euler systems for adjoint Hilbert modular Galois representations}

\date{\today}
\begin{document}
\begin{abstract}

We prove the existence of Euler systems for adjoint modular Galois representations using deformations of
 Galois representations coming from Hilbert modular forms and relate them to $p$-adic $L$-functions under  a conjectural 
 formula for the Fitting ideals of some equivariant congruence modules for abelian base change.
\end{abstract}


\setcounter{section}{0}
\maketitle

\section{Introduction}
After the introduction of the notion of Euler systems  by Kolyvagin  \cite{Kolyvagin} as a powerful tool to understand
 the structure of Selmer groups, a systematic study of then by Perrin-Riou \cite{BPR1} and Rubin \cite{Rubin}
has led naturally to the notion of Euler systems of rank $d$ for some positive integer $d$. In \cite{Rubin-Stark}, Rubin shows that Stark's
 type conjectures give some evidence of the existence of such systems
but so far there was no construction when $d>1$ without using the validity of the  corresponding Main Iwasawa-Greenberg conjecture.

 This work is one of a series of papers (see for example \cite{U-euler-eis}) in which the author is investigating the construction of Euler 
 systems via the study of congruences between automorphic forms of various levels and 
 weights. The main purpose of this note is to give an illustration of this principle for adjoint modular Galois representations and must 
 be seen as an example of a very general construction.  In each situation, some technical difficulties arise and can or cannot be overcome depending of what is known about the 
 structure of certain modules over the Hecke algebras involved. Nevertheless, it is the author's conviction that these principles shed light on the Iwasawa theory
 of the Galois representation at play and will eventually lead to the proof of new Iwasawa Main Conjectures by a successful study of these congruences.
 In the situation of this paper, we will also show that Euler systems
 of rank $d>1$ can be constructed by this technique. Note that despite there are other works giving evidences of the existence of
  Euler system of higher ranks based on the knowledge of the Main conjectures, there is not yet any geometric constructions of such.

To describe the construction done in this paper, let us introduce some notations.
 Let $p$ be an odd prime and $F$ be a totally real number field of degree $d$ over the rationals $\Q$. Let $f$ be a $p$- nearly-ordinary  Hilbert cuspidal eigenform.
Let us  denote by $\rho_f$ the Galois representation attached to $f$: 
$$\rho_f\colon G_F\rightarrow GL_O(T_f)$$
where $T_f$ is a $O$-free module of rank $2$ for some finite extension  $O$ of the ring of $p$-adic integers $\Zp$ and $G_F$ is the absolute Galois group of $F$.
 We will assume throughout the paper that  the residual representation $\bar\rho_f$ is 
absolutely irreducible. This representation is nearly ordinary at each place $v$ dividing $p$, which means for such $v$ there exists
a $O$-direct factor $Fil_v^+T_f$ of rank $1$ which is stable under the action of the decomposition subgroup  $D_v\subset G_F$ at $v$. We will also assume
 that $\rho_f$ is $v$-distinguished\footnote{The trace of the residual representation restricted to $D_v$ is the sum of two characters which are distinct modulo the uniformizer $\varpi$ of $O$.} at each place 
$v$ dividing $p$. We write  $F_v$ for he completion of $F$ at $v$ and $d_v:=[F_v:\Qp]$.

We denote by $ad(\rho_f)\subset End_O(T_f)$ the adjoint representation 
on the endomorphisms of  $T_f$ having trace $0$.
The filtrations  $Fil_v^+T_f$ on $T_f$ induce for each $v$ a three steps  filtration on $ad(\rho_f)$ :
$$\CF^+_v\subset \CF^0_v\subset\CF^-_v=ad(\rho_f)$$
with rank 1 graded pieces. We denote by $Gr^0_v:=\CF^0_v/\CF^+_v$ and fix an isomorphism of $D_v$-module $Gr^0_v\cong O$.
We consider the restriction map at $p$:
$$res_p\colon H^1(F,ad(\rho_f))\rightarrow \bigoplus_{v|p}H^1(F_v,ad(\rho_f))$$
Finally recall that for any Galois representations $V$ of $G_F$ and $S$ a finite set of finite places, we denote by $L^S(V,s)$ the corresponding $L$-function
defined as the Euler product:
$$L^S(V,s):=\prod_{v\notin S} P_v((q_v^{-s};V)^{-1}$$
where $P_v(X,V):=det(1-XFrob_v;V^{I_v})$ with $Frob_v\in D_v$ a geometric Frobenius, $I_v\subset D_v$ the inertia 
subgroup at $v$ and $q_v$ the cardinality of the residue field at $v$. We denote by $\Gamma(V,s)$ the corresponding $\Gamma$-factor.

Let us denote by $\otimes^\Q_F\bar\rho_f$ the tensor induction from $G_F$ to $G_\Q$ of the residual representation $\bar\rho_F$.
The beginning of this work starts by observing that  Hida theory for Hilbert modular forms  can provide a proof of the following
theorem.
\begin{thm}\label{thm1} Let us assume that $\otimes^\Q_F\bar\rho_f$ is absolutely irreducible.
Then, there exists a  canonical element $z_f\in\bigwedge_O^dH^1(F,ad(\rho_f))$  defined up to an element  in $O^\times$
such that 
\begin{itemize}
\item[(i)]  $\wedge ^d res_p (z_f)$ belongs to  $ \bigotimes_{v|p}\bigwedge_O ^{d_v} H^1(F_v,\CF^0_v)$ .
\item[(ii)] The image of $\wedge ^d res_p (z_f)$ by the map
$$\bigotimes_{v|p}\bigwedge_O^{d_v} H^1(F_v,\CF^0_v)\longrightarrow \bigotimes_{v|p}\bigwedge_O^{d_v} H^1(I_v,Gr^0_v)^{D_v/I_v}\cong O$$
is equal to 
$$\xi_f:=\frac{\Gamma(ad(\rho_f),1)L^{S_f}(ad(\rho_f),1)}{\Omega_f^\Sigma\Omega_f^{\Sigma_F\bslh\Sigma}}$$
\end{itemize}
where $S_f$ is the set of finite places containing those where $\rho_f$ is ramified and $(\Omega_f^\Sigma)_{\Sigma\subset\Sigma_F}$ 
are the canonical complex periods attached to the Hilbert modular form $f$ in \cite{Dimitrov-Ihara}.
\end{thm}

The proof of this result follows from  examining the first fundamental exact sequence of K\"ahler differentials of the
universal ordinary Hecke algebra and the interpretation of the latter in terms of Galois cohomology classes
and congruence modules.
In particular, it uses the fact due to the works of  Hida and Wiles \cite{Hi81,Wi95}  in the case $F=\Q$ and
 Dimitrov  \cite{Dimitrov-Ihara} in general  that $\eta_f$ is the size of congruence module
$\wp_f/\wp_f^2$ where $\wp_f=Ker(\lb_f)$ with $\lb_f$ is the homomorphism of the
cuspidal Hecke algebra of same weight and level as $f$ that gives the Hecke eigenvalues associated to $f$.

Before stating the main result of this work which is a generalization of  Theorem \ref{thm1}, let us introduce some more notations. 
 For any number field or $p$-adic field $K$, 
we denote by $H^{cyc}$ the $\Zp$-cyclotomic extension of $K$. For any Galois module $M$ over $G_K=Gal(\bar K/K)$, let $H^1(K,M)$ (resp. $H^1_{Iw}(K,M)$) the Galois cohomology of $M$ 
(resp. the Iwasawa Galois cohomology of $K$ which is defined as the projective limits under the norm maps of 
$H^1(K',M)$ for finite extensions $K'/K$ with $K'\subset K^{cyc}$). We also write $\Lambda:=O[[\Gal(F^{cyc}/F)]]$.

Let us recall that the un-primitive Selmer group $Sel^S(F^{cyc},ad(\rho_f))$ attached to $ad(\rho_f)$ is defined as the kernel of the restriction map:

$$H^1(Gal(F_S/F^{cyc}),ad(\rho_f)\otimes_\Zp \Qp/\Zp)\rightarrow \bigoplus_{v|p}H^1(F^{cyc}_v,ad(\rho_f)/\CF^+_v\otimes_\Zp\Qp/\Zp) $$

\noindent
where $F_S$ is the maximal extension of $F$ unramified away from $S$ and $p$. It is known thanks to the work
 of Wiles, Taylor-Wiles, Fujiwara and others on the modularity of Galois representations (often known ad " $R=T$" Theorems), that this Selmer group is of co-torsion over the Iwasawa algebra $\Lambda$.
We denote by $X^S(F^{cyc},ad(\rho_f))$ its Pontrjagin dual.  Thanks to the work of R. Greenberg \cite[Prop. 4.1.1]{Greenberg}, it s known that
this module does not contain any non trivial finite submodule and therefore its fitting ideal is equal to its characteristic ideal and is therefore principal.
We fix $\CL^{S,alg}_f \in\Lambda$  one of its generators.

By extending the arguments of Theorem 1 for the base change of $f$ to the totally real abelian extensions $E$ of $F$, 
we obtain the following Theorem\footnote{See \S\ref{zeta-para} for  the definition of $\bigwedge$}
 
 \begin{thm}\label{thm2}
There exists a (non trivial) Iwasawa-Euler system of rank $d$ for $ad(\rho_f)$. In other words,
for  totally real fields $E$ running in an $S$-admissible set of  abelian extensions of $F$ which are unramified  above $S$ (see def. \ref  {S-admissible}), there exists an element $Z_{f,E}\in\bigcap^d_{\Lambda[Gal(E/F)]}H^1_{Iw}(E,ad(\rho_f))$  
such that for any extensions $E$ and $E'$ such that $E'\supset E\supset F$, we have

 $$Cores^{E'}_{E}(Z_{f,E'})=\prod_{v\in S(E'/E)}P_v(q_v^{-1}Frob_v,ad(\rho_f)).Z_{f,E}$$
where $S(E'/E)$ are the set of finite places of $F$ that ramifies in $E'$ but not in $E$. Moreover,

\begin{itemize}
\item[(i)] $\wedge^d res_p (Z_{f,F})$ belongs to  $ \dis     \bigotimes_{v|p}\bigcap_{\Lambda}^{d_v} H^1_{Iw}(F_v,\CF^0_v)$ 
\item[(ii)] The image of $\wedge^d res_p (Z_{f,F})$ by the map
$$\bigotimes_{v|p}\bigcap_{\Lambda}^{d_v} H^1_{Iw}(F_v,\CF^0_v)\rightarrow \bigotimes_{v|p}\bigcap_{\Lambda}^{d_v} H^1_{Iw}(I_v,Gr^0_v)^{D_v/I_v}\cong \Lambda $$
is equal to $ \CL^{S,alg}_f $.
\end{itemize}
\end{thm}

The Iwasawa-Zeta elements $Z_{f,E}$ are defined as the projective limit over $n$ of elements $z_{f,E_n}\in \bigwedge^d_{O[\Delta_{E_n}]}H^1(E_n,ad(\rho_f))$; here  $E_n$ denotes  the field of degree $p^n$ over $E$ inside the $\Zp$-cyclotomic extension $E^{cyc}$ of $E$.
In order to construct the elements $z_{f,E} \in\bigwedge^d_{O[\Delta_E]}H^1(E,ad(\rho_f))$  when $E$ runs in a certain set of totally real abelian extensions of
 $F$ with a version of (ii) similar to Theorem 1,
we need to understand the structure of $\wp_{f_E}/\wp_{f_E}^2$ as a module over $O[\Delta_E]$. Here $\Delta_E:=Gal(E/F)$, $f_E$ is the base change to $E$ of the
Hilbert modular form $f$ and $\wp_{f_E}$ is determined similarly as $\wp_f$. The link between the Iwasawa theory of the Selmer 
groups attached to $ad(\rho_f)$, the congruence modules  $\wp_{f_E}/\wp_{f_E}^2$ and deformation theory of the Galois representation $\bar\rho_f$ has been systematically studied by Hida in \cite{HidaSelmer3}. 
In particular, it can be seen from deformation theory that we have a canonical surjection
$$X^S(E^{cyc},ad(\rho_f))\twoheadrightarrow \wp_{f_E}/\wp_{f_E}^2$$
and in the proof of  Theorem \ref {thm2}, we use that the Fitting ideal  of the dual Selmer groups  $X^S(E^{cyc},ad(\rho_f))$ which therefore annihilates    $\wp_{f_E}/\wp_{f_E}^2$ 
 behave well when $E$ varies (see Proposition \ref{selmerprop}).

B. Perrin-Riou \cite[App. B]{BPR1} and K. Rubin \cite{Rubin} have developed some arguments to extract  Euler systems of rank one from one 
of higher rank. However, our method allows us to obtain directly rank one Euler systems with prescribed local conditions at places
dividing $p$. More precisely, we have the following result.
\begin{thm} \label{thm3}For each  place $v|p$ and $h_v\in   H^1(I_{F_v},O)^{G_{F_v}/I_{F_v}}$, then there exists a system of
classes 
 $c_E\in  H^1_{Iw}(E,ad(\rho_f))$
 with $E$ running in an $S$-admissible  set of abelian extensions  of $F$
 such that 
$$Cores^{E'}_{E}(c^{h_v}_{E'})=\prod_{w\in S(E'/E)}P_v(q_v^{-1}Frob_v,ad(\rho_f)).c^{h_v}_{E}$$
Moreover
 \begin{itemize}
\item[(i)]  $ res_w (c^{h_v}_F)$ belongs to  $ H^1_{Iw}(E_w,\CF^0_v)$  if  $w|v$ and is trivial if $ w\not{|} v$.
\item[(ii)] The image of $res_v^0 (c^{h_v}_{F})$ by the map
$$ H^1_{Iw}(F_v,\CF^0_v)\rightarrow  H^1_{Iw}(I_v,Gr^0_v)^{D_v/I_v}\cong \Lambda $$
is equal to 
$\xi\cdot h_v$.
\end{itemize}
\end{thm}

To obtain an Euler system with Iwasawa-Zeta elements that we can relate to  $p$-adic L-functions, we would like to know that
$$\CL^{S}_{f,E}\in Fitt_{O[\Delta_E]}(\wp_{f_E}/\wp_{f_E}^2)$$
where $ \CL^{S}_{f,E}\in K[\Delta_E]$ is the unique  element satisfying for every $\chi\in Hom(\Delta_E,\C^\times)$
\begin{equation}\label{interpole-intro}
\chi(\CL^{S}_{f,E})=\frac{G(\chi)^2\Gamma(ad(\rho_f)\otimes\chi,1)L^{S_E}(ad(\rho_f)\otimes\chi ,1)}{\Omega_f^\Sigma\Omega_f^{\Sigma_F\bslh\Sigma}}
\end{equation}
where $S_E$ is the union of $S$ and the set of finite places that ramify in $E/F$ and $G(\chi)$ is the Gauss sum attached to $\chi$.
We make the following conjecture.
\begin{conj} \label{equivariant-conj} For each totally real abelian extension $E/F$ that ramifies away from $S$ or at $p$, 
we have $\CL^{S}_{f,E}\in Fitt_{O[\Delta_E]}(\wp_{f_E}/\wp_{f_E}^2)$.
\end{conj}

Then we have:
\begin{thm}\label{thm4}
If we assume  Conjecture \ref{equivariant-conj}, the conclusion of Theorem \ref{thm2}  holds with
$\CL^{S,alg}_f$ replaced by the un-primitive  Coates-Schmidt $p$-adic L-function $\CL^{S,an}_f$.
\end{thm}

A weaker form of the conjecture is the following .
\begin{conj} \label{equivariant-conj-weak} For each totally real abelian extension $E/F$ that ramifies away from $S$ or at $p$, 
 $\CL^{S}_{f,E}$ annihilates  $\wp_{f_E}/\wp_{f_E}^2$.
\end{conj}

Then we have:
\begin{thm}\label{thm5}
If we assume  Conjecture \ref{equivariant-conj-weak}, the conclusion of Theorem \ref{thm3}  holds with
$\CL^{S,alg}_f$ replaced by the un-primitive  Coates-Schmidt $p$-adic L-function $\CL^{S,an}_f$.
\end{thm}

Note that either of these conjectures implies that  the elements $\CL^S_{f,E}$ are integral, in otherwords that they belong to $ O[\Delta_E]$. Note that this fact could be proven using the technics of
 \cite{Hida90}. Assume this integrality statement for all extensions $E_n/E$, we define the equivariant Coates-Schmidt $p$-adic L-function
$$\CL^{S,an}_{f,E}:=\varprojlim_{ n}\CL^{S}_{f,E_n}\in \Lambda =\varprojlim_{ n} O[Gal(E_n/E)].$$
where we write $E_n$ for the extension of $E$ of degree $p^n$ inside $E^{cyc}$. For $E=F$, we just write $\CL^{S,an}_{f}\in\Lambda$. Up to a unit in $\Lambda$, $\CL^{S,an}_{f}$ is the
un-primitive Coates-Schmidt $p$-adic L-function.

Using the  Euler system machinery, we can then deduce, assuming conjecture  \ref{equivariant-conj-weak},
 that the following divisibility holds in $\Lambda$:
\begin{equation}
 \CL^{S,alg}_{f}\;  |\;  \CL^{S,an}_{f}
 \end{equation}
On the other hand, by refining the technics of Iwasawa theory and Euler systems to the equivariant setting, 
it is very likely that one can show that the existence of an Euler system
with the bottom class satisfying the correct property is implied by the divisibility above 
 for all the twist $ad(\rho_f)\otimes\chi$. It would then  follow that the existence of the Euler system 
satisfying the correct bottom class condition is equivalent to Conjecture  \ref{equivariant-conj} . We leave the verification of these expectations to the interested  reader.

This note is organized as follows. The Section 2 is devoted to recall facts about Hida families for Hilbert cuspidal Eigenforms, their Galois representations and congruence modules.
We also state the conjecture about their annihilators. The  section 3 contains the construction and the properties of the cocycles and of the compatible systems 
of cohomology classes and zeta elements.

This work  was undertaken when the author was funded by a NSF grant DMS-1407239. Part of it was carried out when the author visited Universitat
 Polit\`ecnica de Catalunya during the Spring 2019. The author would like to thank this institution and Victor Roger for their hospitality and financial support. 
 Finally, the author wants to thank the referee for pointing out inaccuracies and the comments that helped to improve the writing of this paper.

\newpage
\tableofcontents

\noindent
{\bf Notations.}
Throughout this paper $p$ is a fixed odd rational prime. Let $\bQ$ and  $\bQ_p$ be respectively the 
algebraic closures of $\Q$ and $\Q_p$ and let $\C$ be the field of complex numbers.
We fix embeddings $\iota_\infty:\bQ\hookrightarrow \C$ and
$\iota_p: \bQ\hookrightarrow \bQ_p$.
Throughout this note, we implicitly view $\bQ$ as a subfield of $\C$ and $\bQ_p$ via the 
embeddings $\iota_\infty$ and $\iota_p$ and we fix an identification 
$\bQ_p\cong\C$ compatible with these  embeddings.
We denote by $\epsilon_{cyc}$ the $p$-adic cyclotomic character.

\newpage

\section{Congruences modules for Hilbert Modular Forms}

\subsection{Universal ordinary Hecke algebras} Let $F$ be a totally real field and $O_F$ its ring of integers. We denote by $\Sigma_F$ the set of embeddings of $F$ into $\bQ$.
Let $T=T_F$ be the torus defined over $\Z$ by $T_F(A)=(O_F\otimes_\Z A)^\times$. Let $p$ be an odd prime. We write  $T(\Zp)_{tor}$ for the torsion
 subgroup of $T(\Z_p)$ and $\Gamma_F$ for the subgroup of $T(\Z_p)$
such that $T(\Z_p)\cong T(\Zp)_{tor}\times\Gamma_F$. We have $\Gamma_F\cong \Z_p^d$ with $d=[F:\Q]$ is the degree of $F$ over $\Q$. We fix $O$ a
 finite extension of $\Z_p$, we denote $\Lambda_{F}:=O[[\Gamma_F]]$.

Let  $\lb=(\sum_\sigma k_\sigma .\sigma,\sum_\sigma l_\sigma .\sigma,  ) \in\Z[\Sigma_F]^2$ be an arithmetic weight  (which means that $k_\sigma\geq 2$ for all $\sigma\in \Sigma_F$ 
and $w=k_\sigma+2l_\sigma$  is independent of $\sigma$).
We fix $\n$ a non zero  integral ideal of $O_F$ prime to $p$ and let $K^p(\n) \subset GL_2(\widehat\Z^p\otimes O_F)$ be the  
 subgroup of matrices which are congruent to the identity matrix modulo $\n$
and where we have written $\widehat\Z^p$ for the projective limit of $\Z/M\Z$ for the integers $M$ prime to $p$. We also fix $\omega$ be an id\` ele class character of conductor dividing $\n p^\infty$ and infinity type $|\cdot |^w$.

For each positive  integer $n$, we denote by $K_0(p^n)$ the 
subgroup of $GL_2(O_F\otimes\Z_p)$ of matrices 
which are upper triangular modulo $p^n$ and by  $K_1(p^n)$ its subgroup of those such that the diagonal entries are congruent modulo $p^n$. We will identify 
$K_0(p^n)/K_1(p^n)$ with $(O_F/p^nO_F)^\times$
via the map $\left(\begin{array}{cc} a&b\\c&d\end{array}\right)\mapsto a^{-1}d$.
Let $h_{F,\lb}^{ord}(\n,p^n,\omega)$  be the ordinary Hecke algebra 
of lever $K^p(\n)K_1(p^n)$ and weight $\lb$. We then consider the universal Hecke nearly ordinary Hecke algebra of weight $\lb$ and tame level $K^p$ and action of the center given by $\omega$.
$$\mathbf h^{ord}_F= \mathbf h^{ord}_{F,\lb}(\n):=\varprojlim_{ n} h_{F,\lb}^{ord}(\n,p^n)$$
This is the Hecke algebra denote $h(U_\infty^W,O)$ by Hida in \cite[p. 403]{HidaSelmer3} with $W=Cl_F(p^\infty)$.
From the action of  $K_0(p^n)/K_1(p^n)$ on the space of ordinary forms of weight $\lb$ and level $K^p(\n)K_1(p^n)$, we inherit an action of $\Gamma_F$
on $\mathbf h^{ord}_F$ which therefore has the structure of a $\Lambda_F$-algebra.
For any finite order character $\psi$ of $\Gamma_F$, we denote by $P_{\psi}$ the kernel of the map from $\Lambda_F$ into $\bQ_p$ induced by $\psi$.

\begin{thm}[Hida] \label{Hida} $\mathbf h^{ord}(\n)$ is free of finite rank over $\Lambda_F$. Moreover, for any arithmetic characters $\psi_\lb$, there is a canonical isomorphism
$$\mathbf h^{ord}_{F,\lb}(\n) \otimes\Lambda_F/P_{\psi}\cong h_{F,\lb}^{ord}(\n,  p^r,\psi)$$
where $h_{F, \lb}^{ord}(\n, p^r,\psi)$ denotes the Hecke algebra over $O(\psi)$ generated by the Hecke
 operators acting on the space of ordinary Hilbert modular forms of weight $\lb$, level
$K^p(\n)K_1(p^r)$ and with nebentypus restricted to the image of $\Gamma_F$ into $(O_F/p^rO_F)^\times$ given by 
$\psi$ where $r$ is the smallest integer such that $\psi$ factorizes through that image.
\end{thm}
\proof This follows from \cite{HidaControl}. See for example Corollary 5.3 of   \cite{HidaSelmer3} for $W^?=W^\phi$.\qed

Let now $f$ be a $p$-nearly ordinary Hilbert modular form of tame level $K^p(\n)$, unramified central character at $p$ and weight $\lb$ with
 Hecke eigenvalues contained in $O$ via the embedding of $\bQ$ into $\bQp$. It gives rise to an homomorphism:
$\lambda_f : h_{F,\lb}^{ord}(\n , p^r,\psi)\longrightarrow O$ for some character $\psi$ of level $p^r$.
We will denote $\hat\lambda_f$ the composite of $\lambda_f$ with the canonical surjective 
map $\mathbf h^{ord}_{F,\lb}(\n)\rightarrow  h_{\lb}^{ord}(\n, p^r,\psi)$. We therefore have a map:
$$\hat\lambda_f : \mathbf h^{ord}_{F,\lb}(\n)\longrightarrow O$$
Following H. Hida, we can define two congruence modules. Let $B$ be the unique quotient of $h_{\lb}^{ord}(\n, p^r,\psi)$ such that 
$$h_{F,\lb}^{ord}(\n, p^r,\psi)\hookrightarrow O\times B$$
where the first projection map is $\lambda_f$ and the  second $\lambda_B$ is the canonical map induced by the fact $B$ is a quotient of $h_{\lb}^{ord}(\n, p^r,\psi)$.
Let $\eta_f=ker(\lambda_B)$ which naturally imbeds into $O$ via $\lambda_f$ and let $\wp_f=Ker(\lambda_f)$ which imbeds in $B$ via $\lambda_B$. 
The we can define the congruence modules
$$C_0(f):= O/\eta_f\hbox{ and }C_1(f):=\wp_f/\wp_f^2\cong \Omega_{h_{\lb}^{ord}(\n, p^r,\psi)/O}\otimes_{\lambda_f}O.$$
Via the canonical  isomorphism $O/\eta_f\cong B/\wp_f$, we see that $\eta_f$ annihilates $\wp_f/\wp_f^2$ and $\wp_f/\wp_f^2$ is therefore finite since it is finitely generated over $O$ and $\eta_f\neq 0$.

Let $\hat\wp_{f}:=Ker(\hat\lambda_f)$. The following exact sequence will play a fundamental role in this paper.

\begin{lem} \label{FES} With the notation as before, we have a canonical exact sequence.

$$0\rightarrow \Omega_{\Lambda_F/O}\otimes_{\Lambda_F} O\rightarrow \hat\wp_{f}/\hat\wp_{f}^2\rightarrow \wp_f/\wp_f^2\rightarrow 0$$

\end{lem}
\proof
We write the first fundamental exact sequence of Kahler differentials attached to the maps $O\rightarrow \Lambda_F\rightarrow \mathbf h_{F,\lb}^{ord}(\n )$
which we tensor by $O$ through $\hat\lambda_f$:
$$\Omega_{\Lambda_F/O}\otimes_{\Lambda_F}\h^{ord}(\n )\otimes_{\hat\lambda_f} O\stackrel{(1)}{\longrightarrow}\Omega_{\mathbf h_{F,\lb}^{ord}(\n )/O}\otimes_{\hat\lambda_f} O
\longrightarrow \Omega_{h_{F,\lb}^{ord}(\n )/\Lambda_F}\otimes_{\hat\lambda_f} O\rightarrow 0$$
Since $\mathbf h_{F,\lb}^{ord}(\n )\otimes \Lambda/P_{\psi,\lb}\cong h_{F,\lb}^{ord}(\n, p^r,\psi)$, we easily see that we have
$\Omega_{\mathbf h^{ord}(\n )/\Lambda_F}\otimes_{\hat\lambda_f} O   \cong \Omega_{\mathbf h_{F,\lb}^{ord}(\n )/\Lambda_F}\otimes_{\mathbf h^{ord}(\n )} 
 h_{F,\lb}^{ord}(\n, p^r,\psi)   \otimes_{\lambda_f} O \cong \Omega_{h_{F,\lb}^{ord}(\n, p^r,\psi)/O }\otimes_{\lambda_f} O=\wp_f/\wp_f^2$
and is therefore finite . To deduce that the map (1) is injective, since $\Omega_{\mathbf h^{ord}(\n )/\Lambda_F}\otimes_{\hat\lambda_f} O=\Omega_{\Lambda_F/O}\otimes_{\Lambda_F} O\cong O^d$ 
is torsion free. it is sufficient to prove the injectivity after inverting $p$. Since after inverting $p$, the map (1) is surjective, it is therefore sufficient to prove that
$\Omega_{\mathbf h_{F,\lb}^{ord}(\n )/O}\otimes_{\hat\lambda_f} O[1/p]$ is of rank at least $d$. 
This follows from the fact that  $h_{F,\lb}^{ord}(\n )[1/p]$ is of Krull dimension $d$, by  \cite[p. 40 property ff]{HidaSelmer3}, or because it is free of finite rank over $\Lambda_F[1/p]$ thanks to Theorem \ref{Hida}. Note that the fact that $ \wp_f/\wp_f^2$  is finite says exactly that the map $\Spec h_{F,\lb}^{ord}(\n )\to \Spec\Lambda_F$ is \'etale at $\hat\lb_f$.
\qed

\subsection{Base change}\label{basechange}
Let now $E$ be a totally real abelian extension of $F$. The norm map from $O_{E}^\times$ to $O_F^\times$ induces a map from $\Gamma_E$ to $\Gamma_F$ and therefore
from $\Lambda_{E}$ to $\Lambda_F$ inducing the base change transfer map for characters of $p$-power level. It induces an isomorphsim
$$\Lambda_E/I_E\Lambda_E\cong \Lambda_F$$
where $I_E$ is the augmented ideal of the group ring $A_E:=O[Gal(E/F)]$ acting on $\Lambda_E$ via the Galois action of  $Gal(E/F)$ on $O_E^\times$.
 On the other hand, the imbedding $O_F^\times\hookrightarrow O_E^\times$ gives
$\Lambda_E$ the structure of a $\Lambda_F$-algebra.

Let us identify  $\Sigma_E$ with the  set of embeddings of $E$ into $\bQ_p$ via our fix embedding $\iota_p\colon\bQ\hookrightarrow\bQ_p$.
Then we have a non-canonical isomorphism 
$$\Lambda_E\cong O[[T_\sigma , \sigma\in \Sigma_E]]$$
such that the map $\Lambda_E\to O$ induced by the evaluation at the weight of $f_E$ is induced by $T_\sigma\mapsto 0$ for all $\sigma \in \Sigma_E$.
We may choose and fix  such isomorphisms so that the canonical surjective map $\Lambda_E\rightarrow \Lambda_{E'}$
is induced by the canonical surjection $\Sigma_E\twoheadrightarrow \Sigma_{E'}$.  
This clearly induces the  isomorphism
$$\Omega_{\Lambda_E/O}\otimes O\cong \bigoplus_{\sigma\in I_E} O.dT_\sigma$$
where $\sum_\sigma f_\sigma.dT_\sigma$ is mapped to $\sum_\sigma f_\sigma (0,\dots,0) .dT_\sigma$. Moreover the free action of $Gal(E/F)$ on $\Sigma_E$, shows that
$\Omega_{\Lambda_E/O}\otimes O$ is free of rank $d$ over $A_E$.

We now define the base change homomorphism for the ordinary universal Hecke algebras. Let $\n_E:=\n O_E$. By the existence of solvable base change, there is a canonical algebra homomorphism:
$$\pi^E_F(\n,p^n) \colon  h^{ord}_{E}(\n_{E},p^n)[1/p]\rightarrow h^{ord}_F(\n,p^n )[1/p]$$ 
Let $\CH_E^S$ be the restricted tensor product over $O$ of the spherical Hecke algebras over $O$ for $GL_2(E_w)$ for all finite place $w$ above a place of $F$ not in $S$.
The local base change for unramified representations implies we have a canonical algebra homomorphism $\Pi^E_F(S)\colon \CH^S_E\to\CH^S_F$ such that the following diagram commutes
$$
\xymatrix{
\CH^S_E\ar[r]^{\Pi^E_F(S)} \ar[d]& \CH^S_F\ar[d]\\
h^{ord}_{E}(\n_{E},p^n)[1/p]\ar[r] & h^{ord}_F(\n,p^n )[1/p]
}
$$
Since  the Hecke algebras $ h^{ord}_{E}(\n_{E},p^n)$ are torsion free, it follows from the above diagram that the base change homomorphism 
preserve the integrality of Hecke operators and defines an homomorphism
$h^{ord}_{E}(\n_{E},p^n)\rightarrow h^{ord}_F(\n,p^n )$ . After passing tothe projective limit over $n$, we deduce the base change map for the ordinary universal Hecke algebras:
$$\pi^E_F : \mathbf  h^{ord}_{E}(\n_{E})\rightarrow \mathbf h^{ord}_F(\n)$$ 

 becomes  an homomorphism of  $\Lambda_E$-algebras  for  the $\Lambda_E$-algebra structure of $\mathbf h^{ord}_F(\n )$ 
 induced by base change  transfer map for character $\Lambda_E\rightarrow\Lambda_F$ described above.

Let us denote by $f_E$ the ordinary Hilbert modular form for $E$ defined by $\hat\lb_{f_E}=\hat\lb_f\circ\pi_{E,F}$. 
Recall that there is a natural action of the Galois group $Gal(E/F)$ over  $h^{ord}_{E}(\n .O_{E})$ defined by
$U_{\mathfrak p}\to U_{\mathfrak p^\sigma}$ for $\mathfrak p |p$ and $T_{\mathfrak q}\mapsto  T_{\mathfrak q^\sigma}$ for any prime ideal $\mathfrak q$ of $E$ and $\sigma\in Gal(E/F)$.
It induces an action on the characters of this algebra leaving $\lambda_{f_E}$ invariant. It induces therefore an action 
of $Gal(E/F)$ over $C_1(f_E)=\wp_{f_E}/\wp_{f_E}^2$.  

Let $A_E:=O[Gal(E/F)]$. Let us record the  following lemma 

\begin{lem} \label{FESE} With the notation as before, $\wp_{f_E}/\wp_{f_E}^2$ is $O$-torsion and 
we have a canonical exact sequence of $A_E$-modules

$$0\rightarrow \Omega_{\Lambda_E/O}\otimes_{\Lambda_E} O\rightarrow \hat\wp_{f_E}/\hat\wp_{f_E}^2\rightarrow \wp_{f_E}/\wp_{f_E}^2\rightarrow 0$$

\end{lem}
\proof The proof is similar to Lemma  \ref{FES}. The fact that the maps are equivariant for the action of $\A_E$ is clear.\qed

\subsection{Galois representations}
Let $G_F=Gal(\bQ/F)$ be the absolute Galois group of $F$. 
For each finite place $v$, we denote respectively $G_v,I_v$ and $Frob_v\in G_v$ a decomposition subgroup at $v$, its inertia subgroup and a Frobenius lift at $v$.
To each ordinary Hilbert modular form $f$ for $F$ of weight $\lb$ and tame level $\n p^r$ and nebentypus at $p$ given by $\psi$ as in the 
section above, one can associate a nearly ordinary p-adic continuous Galois representation thanks to to the work of many authors  or example see \cite{Wiles} and \cite[p.406-407]{HidaSelmer3}.
$$\rho_f: G_F\longrightarrow GL_2(O)$$
such that
\begin{itemize}
\item[(i)] $\rho_f$ is unramified at finite places of $F$ not dividing $\n p$.
\item[(ii)] For each place $v$ not dividing $\n p$ and arithmetic Frobenius $Frob_v\in G_F$ at $v$, $Trace(\rho_f(Frob_v))=\lb_f(T_v)$
\item[(iii)] For each place $v|p$, there exists $\alpha_{f,v}\in GL_2(O)$ and characters $\psi_v$ and $\psi'_v$ of $G_v$, such that 
$$\rho_f(g)=\alpha_{f,v}
 \left(\begin{array}{cc}  \psi'_v(g) &*\\ 0& \psi_v(g) \end{array}\right)\alpha_{f,v}^{-1}\quad \forall g\in G_v
$$
such that $\psi_v|_{I_v}$ is the  character of $I_v$ corresponding to $\psi|_{O_{F_v}}$ via local Class Field Theory that identifies $\Frob_v$ with  an uniformizer of $O_{F_v}$.
\item[(iv)] $\det\rho_f=\omega \epsilon_{cyc}$.

\end{itemize}
Let us assume the following irreducibility assumption:
\begin{itemize}
\item[{\bf (Irred)}] The reduction $\bar\rho_f$ of $\rho_f$ modulo the maximal ideal of $O$ is absolutely irreducible.
\end{itemize}

Let $sl_2=sl_2(O)$ be the set of $2\times 2$ matrices $A$ with entries in $O$ such that $tr(A)=0$. We consider the adjoint action $ad(\rho_f)$ of $G_F$ on $sl_2(O)$ via $\rho_f$. It is defined by
$$ad(\rho_f)(g).A=\rho_f(g)A\rho_f(g)^{-1}$$
Note that under the condition{\bf (Irred)}, we have
\begin{equation}\label{H0}
H^0(E,ad(\bar\rho_f))=0
\end{equation}
for any totally real field extension $E/F$.
We will further consider the conditions:

\begin{itemize}
\item[{\bf (Big)}] $ad(\bar\rho_f)$  is absolutely irreducible.
\end{itemize}

\begin{itemize}
\item[{\bf (Dist)}] For each $v|p$, he characters $\psi_v$ and $\psi'_v$ are distinct modulo $p$.
\end{itemize}

\begin{itemize}
\item[{\bf (Indec)}] For each $v|p$, $\rho_f|_{G_v}$ is indecomposable.
\end{itemize}

\subsection{Twisted adjoint L-values}\label{twisted-L}
We are interested in the nature of the fitting ideal of $\wp_{f_E}/\wp_{f_E}^2$ over $A_{E}$. We formulate
 a conjecture about it below. We first introduce some notations.
For any finite order character $\chi$ of $G_F$ and finite set of primes $S$, we  consider  the $S$-primitive  twisted $L$-function
$L^S(ad(\rho_f)\otimes\chi,s)$. It is known to be holomorphic on the whole complexe plane except maybe if $\rho_f$ is dihedral, a case we exclude by our hypothesis.
For any subset $\Sigma\subset\Sigma_F$, M. Dimitrov has introduced in \cite{Dimitrov-Ihara} some canonical periods $\Omega_f^\Sigma\in\C^\times/(O\cap\bQ)^\times $ so that
$$\CL^{S}(ad(\rho_f)\ot\chi ,1):=
G(\chi )^2\frac{\Gamma(ad(\rho_f)\ot\chi,1)L^S(ad(\rho_f)\ot\chi , 1)}{\Omega_f^\Sigma\Omega_f^{\Sigma_F\bslh\Sigma}  } \in K(\chi )$$
where $K(\chi )\subset\bQp$ is the extension of $K$ generated by the values of $\chi $ and $G(\chi )$ is the Gauss sum attached to $\chi $.

Let $S_E$ be the set of places of $F$ dividing $\n p$ and those ramifying in the extension $E/F$. For each character
 $\chi \colon Gal(E/F)\to\bQ^\times\subset\bQp^\times$, let $e_\chi \in K(\chi ) \otimes_O A_E$ be the idempotent projecting any $K[Gal(E/F)]$-modules
on its $\chi $-isotypical component. Let $A_E:=O[Gal(E/F)]$. We then define
$$\CL_E^{an}(ad(\rho_f)):=\sum_{\chi }
\CL^{S_E}(ad(\rho_f)\ot\chi ,1).e_\chi  \qquad \in  \bQ_p\otimes_O A_E$$
where the sum runs over all the characters of $Gal(E/F)$. The part (iv) of the following conjecture is an equivariant formulation 
of the congruence number formula established by Hida in the eighties.

\begin{conj} \label{equiv-conj}
For any totally real abelian extension $E/F$ that ramifies away from $S$ or at $p$, we have:

\begin{itemize}
\item[(i)] $\CL_E^{an}(ad(\rho_f))\in A_E$,
\item[(ii)]  $\CL_E^{an}(ad(\rho_f))$ annihilates $\wp_{f_E}/\wp_{f_E}^2$,
\item[(iii)] $\CL_E^{an}(ad(\rho_f))$ belongs to the Fitting ideal of $\wp_{f_E}/\wp_{f_E}^2$ over $A_E$.
\item[(iv)] $\CL_E^{an}(ad(\rho_f))$  generates  the Fitting ideal of $\wp_{f_E}/\wp_{f_E}^2$ over $A_E$ 

\end{itemize}
\end{conj}

The part (i) of the above conjecture should not be difficult to prove using the  integral representation of these twisted $L$-values using Eisenstein series and Theta series
of half integral weight as it was done by Hida in \cite{Hida90} for $F=\Q$. We can think of (ii) as an automorphic Stickelberger Theorem and we will see that
if we assume (ii) for many fields $E$, we can construct Euler system of rank one for $ad(\rho_f)$. However part (iii) which implies part (ii) is stronger and allows to construct an
  Euler systems for $ad(\rho_f)$ of rank $d$. When $E/F$ is quadratic, J. Tilouine and the author have proved this conjecture in \cite{TU} which was formulated
in the quadratic case by H. Hida in a slightly different way.

\section{Congruences and  Euler systems}

\subsection{Big Hecke rings and Galois representations}
Let $E/F$  like in section \ref{basechange}.  
Let $\mathbf T_E$ be the local component of $\mathbf h^{ord}_E(\n)$ associated to the maximal ideal containing $Ker(\hat\lb_{f_E})$. 
It is flat over $\Lambda_E$ and it is equipped with a natural action of  $ Gal(E/F)$ such that the natural map
$\Lambda_E\to \mathbf T_E$ is $Gal(E/F)$-equivariant.

Then, there exists a Galois representation
$$\rho_{\mathbf T_E}: G_E\longrightarrow GL_2(\mathbf T_E)$$
such that
\begin{itemize}
\item[(i)] $\rho_{\mathbf T_E}$ is unramified at finite places of $F$ not dividing $\n p$.
\item[(ii)] For each place $v$ not dividing $\n p$ and arithmetic Frobenius $\Frob_v\in G_E$ at $v$, $Trace(\rho_{\mathbf T_E}(\Frob_v))=T_v$
\item[(iii)]For each place $v|p$ of $F$ and each place $w|v$ of $E$, there exists $\alpha_w\in GL_2(\mathbf T_E)$ and characters $\Psi_w$ and $\Psi'_w$ of $G_v$, such that 
$$\rho_f(g)=\alpha_{w}
 \left(\begin{array}{cc}  *&*\\  &\Psi_w(g) \end{array}\right)\alpha_w^{-1}\quad \forall g\in G_w
$$
and  $\Psi_w|_{I_w}=\psi_v|_{I_w} \kappa_w$, where $\kappa_w$ is the  character of $I_w\stackrel{Art}{\longrightarrow} O_w^\times  \hookrightarrow (O_E\otimes\Zp)^\times  \twoheadrightarrow \Gamma_E \subset\Lambda_E^\times$ where $Art$ is the the restriction to $I_w$ of the isomorphsim of  local Class Field Theory $G_w^{ab}\cong E_w^\times$ that identifies $\Frob_w$ with  an uniformizer of $O_{E_w}$.

\item[(iv)] $\det\rho_{\T_E}=\omega \epsilon_{cyc}$.
\end{itemize}
\subsection{Construction of cocycles}\label{cocycle-cons}
Let $\hat\wp_{f_E}=Ker(\hat\lb_{f_E} : \mathbf T\rightarrow O)$. 
We denote  $\q_E:=\wp_{f_E}/\wp_{f_E}^2$ and $\hat\q_E:=\hat\wp_{f_E}/\hat\wp_{f_E}^2$ , and we see them 
as Galois modules with action of $G_F$ given via the natural action of $Gal(E/F)$ on $\mathbf T_E$.
Note that the exact sequence of Lemma \ref{FESE}, induces an exact sequence of $O[Gal(E/F)]$-modules:
\begin{equation}\label{FESEequiv}
0\to \Omega_{\Lambda_E/O}\otimes O\to \hat\q_E\to\q_E\to 0
\end{equation}
where $\Omega_{\Lambda_E/O}\otimes O$ is free of rank $d$ over $A_E=O[Gal(E/F)]$.

 Let us denote by $sl_2(O)$ the Lie algebra over $O$ of $2\times 2$ matrices of trace $0$.
We will most of the time  use the notation $ad(\rho_f)$ instead of $sl_2(O)$ to emphasis that it is a Galois module for the adjoint action of $\rho_f$ given by
$$ g.X:=\rho_f(g)X\rho_f(g)^{-1}\qquad \forall g\in G_F,\quad \forall X\in sl_2(O)$$
For $g\in G_E$, we define an element $x_{f,E}(g):= x(g)\in sl_2(O)\otimes_O\hat\wp_{f_E}/\hat\wp_{f_E}^2$ by the equality
$$x(g):=\rho_{\mathbf T}(g)\rho_f(g)^{-1}  -\mathbf 1_2 \pmod{\hat\wp_{f_E}^2}$$

\begin{lem}
The map $g\mapsto x_{f,E}(g)$ defines a cocycle in $Z^1(G_E, ad(\rho_f)\otimes\hat\q_E)$.
Moreover, for any $\tau\in Gal(E/F)$, we have
$$x^{\tau}(g)=\tau^{-1}.x(g)\; ,\quad\forall g\in G_E$$
where on the left hand side the action of $\tau$ is the natural one on the set of cocycles and 
the right hand side action is the one on $ \hat\q_E $.
\end{lem}
\proof 
From the definition of $x(g)$ we have
$$\rho_{\mathbf T}(g)=(\mathbf 1_2+x(g))\rho_f(g) \pmod{ \hat\wp_{f_E}^2} $$ 
Since $\rho_\T$ and $\rho_f$ have same determinant, we have  $det((1+x(g))=1 \pmod{ \hat\wp_{f_E}^2}$ and therefore $x(g)$ has trace zero modulo $\hat\wp_{f_E}^2$.
Now, an easy computation provides that
\begin{eqnarray*}
x(gg')&=&x(g)+\rho_f(g)x(g')\rho_f(g)^{-1}+x(g)\rho_f(g)x(g')\rho_f(g)^{-1}\\
& =&x(g)+\rho_f(g)x(g')\rho_f(g)^{-1} \pmod{ \hat\wp_{f_E}^2}
\end{eqnarray*}
which implies the first claim.

Let now $\tau\in Gal(E/F)$. By definition of he action of $\tau$ on $\mathbf T$, we have the matrices identity
\begin{eqnarray}
\label{action1} \tau^{-1} (\rho_{\mathbf T}(g))=A_\tau \rho_{\mathbf T}(\tau ^{-1} g\tau)A_\tau ^{-1}
\end{eqnarray}
for some matrix $A_\tau \in GL_2(\mathbf T) $. Since $\tau.\lb_{f_E}=\lb_{f_E}$, from the 
equality $\rho_f (\tau g\tau^{-1})=\rho_f(\tau)\rho_f(g)\rho_f(\tau)^{-1}$ and
the fact $\rho_f$ is residually irreducible, we can easily see that we may choose $A_\tau$ so that
$$A_\tau\equiv \rho_f(\tau) \pmod{\hat\wp_{f_E}}$$
Therefore, since $\tau^{-1}.\rho_f(g)=\rho_f(g)$ for all $g\in G_E$, we have the following identities modulo $\hat\wp_{f_E}^2$
\begin{eqnarray*}
(\mathbf 1_2+\tau^{-1}.x(g))\rho_f(g)\tau^{-1}.  & \equiv &  
 (\mathbf 1_2+A_\tau x(\tau^{-1} g\tau  ) A_\tau  ^{-1} )A_\tau \rho_f(\tau  ^{-1} g\tau )A_\tau ^{-1}  \\
&\equiv &   (\mathbf 1_2+\rho_f(\tau)  x(\tau ^{-1} g\tau ) \rho_f(\tau) ^{-1})\rho_f(g) 
\end{eqnarray*}
Therefore we deduce:
$$\tau^{-1}.x(g)\equiv \rho_f(\tau)  x(\tau ^{-1} g\tau ) \rho_f(\tau) ^{-1} \pmod{\hat\wp_{f_E}^2}$$
which proves the second claim.
\qed

\begin{coro}
The cocycle $x_{f_E}$ descends canonically to a cohomology class 
$$\tilde x_E\in H^1(F, ad(\rho_f)\otimes_O\hat\q_E).$$
\end{coro}
\proof It is clear that $x_{f_E}$ defines a cohomology class in $H^1(E, ad(\rho_f)\otimes_O\hat\q_E)$. 
This class is invariant by $Gal(E/F)$ by the previous Lemma. The 
 inflation restriction exact sequence  gives
 \begin{eqnarray*}
 0\to H^1(Gal(E/F), H^0(E,ad(\rho_f))\otimes_O\hat\q_E)  \to H^1(F, ad(\rho_f)\otimes_O\hat\q_E) )\to  \\
  \to H^1(E, ad(\rho_f)\otimes_O\hat\q_E) )^{Gal(E/F)}\to  H^2(Gal(E/F), H^0(E,ad(\rho_f)\otimes_O\hat\q_E))
  \end{eqnarray*}
This implies that $H^1(F, ad(\rho_f)\otimes_O\hat\q_E) )\cong 
  H^1(E, ad(\rho_f)\otimes_O\hat\q_E) )^{Gal(E/F)}$, and that our
result follows since it is easily seen that $H^0(E,ad(\rho_f))\otimes_O\hat\q_E) =0$ from the irreducibility of  residual representation 
 $ad(\bar\rho_f)|_{G_E}$.\qed

In what follows, we will consider the map
$$\mathbf c_E: \hat\q_E^*=Hom_{A_E}(\hat\q_E,A_E)\longrightarrow  H^1(F,ad(\rho_f)\otimes \Lambda_E)=H^1(E,ad(\rho_f))$$
sending $\phi\in Hom_{A_E}(\hat\q_E,A_E)$ to the cohomology class of the cocycle 
$$g\mapsto (id_{sl_2(O)}\otimes \phi)\circ \tilde x_{E}(g).$$

\subsection{Local non-triviality}

Let $L$ be a finite extension of $\Qp$ with absolute Galois group $G_L$. 
We denote by $I_L\subset G_L$ the inertia subgroup of $G_L$ and by $\Frob_L \in G_L$ a Frobenius element. Let
$ S_L$ be the universal deformation ring representing the deformation functor $\CF_1$ of the trivial character 
of $G_L$ defined on the category of local noetherian $O$-algebra. Let $\kappa_L$ be the corresponding 
 universal deformation
 $\kappa_L\colon G_L\to S_L$. We also denote by $\Gamma_L$  the $p$-torsion free part of $O_L^\times$.
 The following lemma is well known.
\begin{lem}  \label{local-def}With the previous notations, the following holds.
\begin{itemize}
\item[(i)] For any $O$-deformation $\chi\in\CF_1(O)$, we have a canonical isomorphism  $S_L\cong O[[T]][[\Gamma_L]]$  with $1+T=\kappa_L(\Frob_L)$ 
such that the map induced by reducing modulo the augmentation ideal of $\Gamma_L$ and $T$ to $0$ corresponds to the character $\chi$. 
\item[(ii)]  $\kappa_L|_{I_L}$ takes values in $O[[\Gamma_L]]^\times$ and is induced by the Artin reciprocity map $Art_L\colon G_L^{ab}\cong L^\times$.
\item[(iii)] For any $O$-deformation $\chi\in\CF_1(O)$, we have canonical isomorphisms:
$$\CF_1(O[\epsilon](\epsilon^2))=Hom_{S_L}(\Omega_{S_L/O}, O)\cong H^1(L,O)$$
where $O$ is seen as a $S_L$-module via the map $S_L\to O$ induced by $\chi$. In particular,
it is  free of rank $[L:\Qp] +1$. Moreover, we have  the following 
 commutative square.
$$ \xymatrix{Hom_{S_L}(\Omega_{S_L/O}, O)\ar[d]\ar[r]^{\qquad\sim}  & H^1(L,O)\ar[d]\\
 Hom_{O[[\Gamma_L]]}(\Omega_{O[[\Gamma_L]]/O}, O)\quad \ar[r]_{\qquad\delta_L}^{\qquad\sim} &\quad  H^1(I_L,O)^{G_L/I_L}
 }$$
 where the horizontal arrows are isomorphism and where the left arrow is induced by the inclusion $O[[\Gamma_L]]\subset S_L$
 and $O$ is an $S_L$-module by the character $\chi$.
 \item[(iv)] These isomorphisms are norm compatible for any extension $L'/L$.
\end{itemize}
\end{lem}
\proof This is an elementary exercise in deformation theory using the reciprocity law of Local Class Field  theory.
\qed


Let $v$ be a place of $F$ dividing $p$ and let $F_v$ be the corresponding 
completion. Let $ad(\rho_f)=\CF_v^-\supset \CF^0_v\supset \CF^+_v $ be the $G_{F_v}$ 
stable filtration with graded pieces $Gr^-_v, Gr^0_v$ and $Gr^+_v$ of rank $1$ over $O$ where
$$\CF_v^0=\{ A\in sl_2(O) | A=\alpha_v^{-1} 
 \left(\begin{array}{cc} a & b \\0& -a \end{array}\right) \alpha_v, \;\hbox{ for  some }a,b\in O\}.$$
 and
 $$\CF_v^+(ad(\rho_v))=\{ A\in sl_2(O) | A=\alpha_v^{-1} 
 \left(\begin{array}{cc} 0 & b \\0 & 0 \end{array}\right) \alpha_v, \;\hbox{ for  some }b\in O\}.$$
 
 On the other hand, we have an isomorphism
\begin{equation}\label{decomp}
\Omega_{\Lambda_E/O}\otimes O \cong  \bigoplus_{w|p}\Omega_{O[[\Gamma_{E_w}]]/O}\otimes O
\end{equation} 
For any $\phi\in Hom_{A_E}(\hat\q_E,A_E)$, we denote by $\phi_w$, the 
restriction of $\phi$ to $\Omega_{O[[[\Gamma_{E_w}]]/O}\otimes O\subset \Omega_{\Lambda_E/O}\otimes O\hookrightarrow\hat\q_E$.

\begin{lem}\label{local-nontriv0}
For each $\phi\in Hom_{A_E}(\hat\q_E,A_E)$ and $w$ a place of $E$ above $v$, we have

\begin{itemize}
\item[(i)]  We have $\mathbf c_w(\phi):=\mathbf c_E(\phi)|_{G_{E_w}}\in H^1(E_w,\CF^0_v)$.
 
\item[(ii)]   Let $\mathbf c^0_w(\phi)$ be the image of $\mathbf {c}_w(\phi)$
in $H^1(I_w,Gr^0_v)\cong H^1(E_w,O)$. Then 
$$\mathbf c^0_w(\phi)=\delta_{E_w}(\phi_w)$$
where $\delta_{E_w}$ is the isomorphism of Lemma \ref{local-def}.(iii)
 \end{itemize}
\end{lem}
\proof
Since $\rho_{\T_E}|_{G_w}$ is nearly ordinary, for $g\in G_w$, we have
$$
\rho_{\T_E}(g) =  \alpha_{w} \left(\begin{array}{cc} * & * \\& \Psi_w(g) \end{array}\right)\alpha_{w} ^{-1}
$$
where thanks to {\bf (Indec)}, we may choose $\alpha_w$ so that $\alpha_{f,v}=\lb_{f_E}(\alpha_w)$ with
$$\rho_f(g)=\alpha_{f,v}
 \left(\begin{array}{cc}* &*\\ 0& \psi_v(g) \end{array}\right)\alpha_{f,v}^{-1}\quad \forall g\in I_v
$$

Let $X_w \in sl_2(O)\otimes_O\hat\wp_E$ such that $X_w\equiv \beta_w^{-1}\alpha_v-\mathbf 1_2\pmod{\hat\wp_E^2}$
For all $g\in G_w$, then an elementary computation gives the equality in $sl_2(O)\otimes\hat\wp_E/\hat\wp_E^2$:
\begin{eqnarray*}
\alpha_w^{-1}x(g)\alpha_w \equiv \left(\begin{array}{cc} * &*\\ 0&\Psi_w\psi_v^{-1}(g) \end{array}\right)-\mathbf 1_2
+ \alpha_{f,v}^{-1}[ad(\rho_f(g))(X_v)-X_v] \alpha_{f,v} \pmod {\hat\wp_E^2}
\end{eqnarray*}

The cocycle defined by the  second term in the above equation being a co-boundary, we may assume that $X_v=0$ in the computation.
We see that from this expression and the definition of the cocycle that (i) follows. Now we look at the restriction to the inertia subgroup $I_w\subset G_w$.
Recall that $\Psi_w|_{I_w}$ is the universal deformation of the character $\psi|_{I_w}$, where we see $\psi$ as a character of $I_v$ for $v$ the place of $F$ below $w$ associated to $\psi$ via the Reciprocity Law of Local  Class Field Theory. In other words, $\Psi_w|_{I_w}=\psi_v|_{I_w}\kappa_{E_w}$.
Now, for $\phi_w\in Hom_O(\Omega_{O[[\Gamma_{E_w}]]},O)$ , we have
$$\delta_{E_w}(\phi_w)(g)=\phi_w(\kappa_{E_w}(g)-1)$$
where we see $\phi_w$ has a linear form on the augmentation ideal $I_{\Gamma_{E_w}}$ of $O[[\Gamma_{E_w}]]$ via
$I_{\Gamma_{E_w}}\twoheadrightarrow I_{\Gamma_{E_w}}/I_{\Gamma_{E_w}}^2\cong\Omega_{O[[\Gamma_{E_w}]]/O}$.
Our claim now follows easily  from interpreting the last three equations in the definition of $\mathbf c_w^0(\phi)$
after identifying $Gr^0_v$ with the trivial $G_{F_v}$-module via the map
$$\alpha_{f,v}.   \left(\begin{array}{cc}-a &*\\ 0& a \end{array}\right)  \alpha_{f,v}^{-1}\mapsto a.$$
\qed
The following corollary is immediate.
\begin{coro}\label{local-nontriv} The map $ \oplus_{w|p}Res_w\circ\mathbf c_E$ induces an isomorphism of $A_E$-module:
$$ \bigoplus_{w|p}Hom_O(\Omega_{O[[\Gamma_{E_w}]]/O}, O)\cong \bigoplus_{w|p} H^1(I_{E_w},O)^{G_{E_w}/I_{E_w}}$$
In particular, the map $\mathbf c_E$ is injective. 
\end{coro}

\subsection{Zeta elements}\label{zeta-para}
We construct zeta elements in Rubin's lattice attached to the Galois cohomology of $ad(\rho_f)$. We first recall some definitions.
 Let $A$ be a finite free $O$-algebra. Let $n$ be a positive integer. If $L$ is a $A$-module over a ring $A$, we put following Rubin:
$$\bigcap_A^n L:= (\bigwedge^n_{A}L^*)^*$$ 
where for any $A$-module $M$, we set  $M^*:=Hom_{A}(M,A)$. $L\mapsto \bigcap_A^n L$ defines clearly a covariant functor 
from the category of $A$-modules to itself.
In what follows we take $A=A_E=O[Gal(E/F)]$. Now we consider  the map we have defined in the end of section \ref{cocycle-cons} : 
$$\mathbf c_E: \hat\q_E^*=Hom_{A_E}(\hat\q_E,A_E)\longrightarrow  H^1(F,ad(\rho_f)\otimes \Lambda_E)=H^1(E,ad(\rho_f))$$
and apply the co-variant functor $\cap^d$ :
\begin{eqnarray} \label{zeta-map}
\cap^d\mathbf c_E \colon \bigcap^d_{A_E} \hat\q_E^*\longrightarrow \bigcap_{A_E}^d H^1(E,ad(\rho_f))
\end{eqnarray}
Since $\Omega_{\Lambda_E/O}\otimes O$ is free of rank $d$ over $A_E$, we can fix an isomorphism 
\begin{eqnarray} \label{choice1}
\bigcap^d_{A_E}(\Omega_{\Lambda_E/O}\otimes O)^*\cong A_E
\end{eqnarray}
From the previous corollary, we get easily the following
\begin{coro}\label{wedge-loc-nontriv} We have a commutative diagram
\begin{eqnarray*}
\xymatrix{
\cap^d_{A_E}\q_E^*  \ar[rr]\ar[d]_{\cap^q\mathbf c_E} && \cap^d_{A_E}  (\Omega_{\Lambda_E/O }\otimes O)^* \ar[d]^{ \cong }\cong A_E \\
\cap^d_{A_E} H^1(E, ad(\rho_f)) \ar[rr]^{\cap^q Res_p\quad}&& 
\displaystyle \cap^d_{A_E}. 
H^1(I_{E_w},O)^{G_{E_w}/I_{E_w}} \
}
\end{eqnarray*}
\end{coro}
\proof We just apply the covariant functor $\cap^q_{A_E}$ to the isomorphism of Corollary \ref{local-nontriv}.\qed
We define the local Zeta elements  $z_{E,loc}\in \cap^d_{A_E} \oplus_{w|p} H^1(E_w,ad(\rho_f))$
the element corresponding to $1$ via the isomorphism of Corollary \ref{wedge-loc-nontriv}
and \eqref{choice1}.
The following lemma will be useful to construct elements in $\bigcap_{A_E}^d H^1(E,ad(\rho_f))
$.
\begin{lem}\label{annihilator}
Let $\tilde\q_E$ be  the cokernel of the map $\bigcap^d_{A_E} \hat\q_E^*\longrightarrow  A_E$
induced by the exact sequence \eqref{FESEequiv} and the isomorphism \eqref{choice1}. Then we have 
$$ Fitt_{A_E}(\tilde\q_E) \supset Fitt_{A_E}(\wp_{f_E}/\wp_{f_E}^2).$$
\end{lem}
\proof
Notice that for any $A_E$-module $M$, we have a canonical isomorphism $M^*=Hom_{A_E}(M,A_E)\cong Hom_O(M,O)=:M'$. Moreover if
we have an injective homomorphism $f\colon X\rightarrow Y$ of $A_E$--modules 
which are $O$-free of same rank over $O$, inducing by duality an injective homomorphism 
$f'\colon Y'\rightarrow X'$. Then $Coker f'$ and $Coker f$ are Pontryagin dual of each other. Moreover if both $X$ and $Y$ are equipped with 
a $O$-linera action of a group $G$ and if $f$ is itself $G$-equivariant, then  the corresponding duality isomorphism is $G$-equivariant.
Which implies that $Fitt_{A_E}(Coker f')=Fitt_{A_E}(Coker f)$.
We consider now the map $f_E\colon \Omega_{\Lambda_E/O}\ot O\rightarrow \hat\q_E$. We have 
$f_E^{**}\colon \Omega_{\Lambda_E/O}\ot O\rightarrow (\hat\q_E)_f$ where $(\hat\q_E)_f$ is the maximal $O$-free quotient of $\hat\q_E$,
 and therefore
$$ Fitt_{A_E}(Coker f_E^{**})\supset  Fitt_{A_E}(Coker f_E)= Fitt_{A_E}(\q_E)=Fitt_{A_E}(\wp_{f_E}/\wp_{f_E}^2)$$
Now by the previous discussion for $f=\wedge^df_E^{**}$, we have
$$Fitt_{A_E}(\tilde\q_E)=Fitt_{A_E}(Coker(A_E\to \wedge^d \q_E^{**}))$$
On the other hand,
$$Fitt_{A_E}(Coker(A_E\to \wedge^d \q_E^{**}))= Fitt_{A_E}(Coker f_E^{**}).$$
The equalities and inclusion above implied our claim.
\qed
We can now give the definition of Zeta elements.
\begin{defi}\label{zeta-def}
For any element $\xi$ annihilating $\tilde\q_E$,  we consider $z_{E,\xi}\in \bigcap_{A_E}^d H^1(E,ad(\rho_f))$ to be the image of
$\xi.1_{A_E}$ view as an element of $\bigcap^d_{A_E} \hat\q_E^*$ by the map $\cap^q\mathbf c_E$ of \eqref{zeta-map}.
\end{defi}
It follows immeditely from Corollaryt \ref{wedge-loc-nontriv}, that
\begin{equation}\label{wedge-localization}
\cap^d Res_p(z_{E,\xi})=\xi. z_{E,loc}
\end{equation} 
\subsection{Proof of Theorem \ref{thm1}}
Notice first that when $E=F$, we have $A_E=O$ and $\bigcap_{A_E}^d=\bigwedge^d_O$.
Let $\eta_f$ be the $O$Fitting idesl of the congruence module attached to $f$. Under our hypothesis, it is known that $Fitt_O(\wp_f/\wp_f^2)=\eta_f$ since
 in that case that the Hecke algebra is complete intersection thanks to the works of Taylor-Wiles-Fujiwara-Dimitrov  (see  \cite{Dimitrov-Ihara} for the most complete statement). Since the cohomology of the Hilbert-Blumenthal variety is free over the Hecke algebra again thanks to the work of Dimitrov  \cite{Dimitrov-Ihara},
 we have  $\eta_f$ is the ideal of $O$ generated by
$$\xi_f:= \frac{\Gamma(ad(\rho_f),1)L^{S_f}(ad(\rho_f),1)}{\Omega_f^\Sigma\Omega_f^{\Sigma_F\bslh\Sigma}}\in O$$
with the notation of Theorem \ref{thm1}.

Now, we just apply  Lemma \ref{annihilator} to the case $E=F$ and the construction of Definition \ref{zeta-def} with $\xi:=\xi_f$.
This gives the construction of the element 
$z_f\in \wedge^d_O H^1(F, ad(\rho_f))$. The point (i) now follows from 
Lemma \ref{local-nontriv0}.(i), and the point (ii) is a direct consequence of the construction and of the relation
\eqref{wedge-localization}.
 \qed

\subsection{Compatible systems for the norm map}
\subsubsection{Construction}
We start by the following lemma

\begin{lem}

Let $E'/F$ be an extension of $F$ contained in $E$. Then the base change map $\T_E\to\T_{E'}$ is surjective. In particular,  the induced map $\hat\q_E\to\hat\q_{E'}$ is surjective and for any $\phi \in \q_{E}^*$,  there exists a unique $\phi'\in \q_{E'}^*$ such that
the following diagram of $A_E$-module commutes
$$
\xymatrix{
\Omega_{\Lambda_{E/O}}\otimes O\ar[r]\ar@{>>}[d]&\Omega_{\mathbf T_E/O}\otimes_{\lb_{f_E}} O\ar@{=}[r]\ar@{>>}[d] & \hat\q_{E}\ar[r]^{\phi}\ar[d]\ar@{>>}[d]&A_E\ar@{>>}[d]^{\pi^{E}_{E'}}\\
\Omega_{\Lambda_{E'/O}}\otimes O\ar[r]\ &\Omega_{\mathbf T_{E'}/O}\otimes_{\lb_{f_{E'}}} O\ar@{=}[r] &\hat\q_{E'}\ar@{-->}[r]^{\phi'}&A_{E'}
}
$$
where $\pi^E_{E'}$ is induced by $Gal(E/F)\rightarrow  Gal(E'/F)$ and the left vertical arrows are induced by the base change maps $\Lambda_E\to\Lambda_{E'}$ and 
 $\mathbf T_E\rightarrow\mathbf T_{E'}$. 
 \end{lem}
 \proof
 Let $R_E$ be the universal nearly ordinary (at places dividing $p$)  deformation ring with fixed determinant equal to $\omega \epsilon_{cyc}$ of $G_E$ unramified at places not dividing  $\n p$ (see \S\ref{kahler-diff} for the precise definition) . It follows from \cite[Prpp. 3.1]{HidaSelmer3}, that $R_E\to R_{E'}$ is surjective.
 Since  for every extension $E$, the canonical map $R_E^\phi\to\T_E$ is surjective. We deduce that $\T_E\to\T_{E'}$ is surjective.  This implies immediately that 
$\hat\q_E\to\hat\q_{E'}$ is surjective. On the other hand, we see easily from the exact sequence \eqref{FESEequiv} for $E$ and $E'$, that the kernel of the surjective map $\hat\q_E\otimes_{A_E} A_{E'}\to\hat\q_E$ is $O$-torsion. We deduce easily the rest of the Lemma from the surjectivity and that last observation.\qed
By the previous lemma, $\phi\mapsto\phi'$ defines a $A_E$-linear map
$$\hat\pi^E_{E'}\colon \hat \q_E^*\rightarrow \hat \q_{E'}^*$$ where the action of $A_E$ onto $\q_{E'}$ is given via the ring homomorphism $\pi^E_{E'}$.
Similarly, we have a map $(\Omega_{\Lambda_E/O}\otimes O)^*\to (\Omega_{\Lambda_{E'}/O}\otimes O)^*$ that we denote $\hat\pi^E_{E'}$ too. 
 We deduce that we have a canonical diagram
\begin{eqnarray}\label{Norm-diag}
\xymatrix{
(\Omega_{\Lambda_E/O}\otimes O)^*\ar[d]^{\hat\pi^E_{E'}}&\ar[l]\hat\q_{E}^*\ar[r]^{\mathbf c_E\qquad\qquad}\ar[d]^{\hat\pi^E_{E'}}\ar[d]&H^1(F,ad(\rho_f)
\otimes A_E)\ar[d]^{id_{ad(\rho_f)}\ot \pi^{E}_{E'}}  \ar@{=}[r]& H^1(E',ad(\rho_f))   \ar[d]^{N^E_{E'}}\\
(\Omega_{\Lambda_{E'}/O}\otimes O)^*&\ar[l]\ \hat\q_{E'}^*\ar[r]^{\mathbf c_{E'}\qquad\qquad}&H^1(F, ad(\rho_f)\otimes A_{E'})\ar@{=}[r]&H^1(E',ad(\rho_f))
}
\end{eqnarray}
where $N^E_{E'}$ stands for the natural co-restriction or Norm  map. We can now choose the isomorphisms \eqref{choice1} so that for any extension $E/E'$ as before,  the following diagram commutes.
$$
\xymatrix{
\bigcap^d_{A_E}\ar[d]\ar[r]^{\qquad\sim}  \Omega_{\Lambda_{E/O}}\otimes O&A_E\ar[d]^{\pi^E_{E'}}\\
\bigcap^d_{A_E}\ar[r]^{\qquad\sim} \Omega_{\Lambda_{E/O}}\otimes O&A_{E'}
}
$$
The following lemma results from the above discussion.
\begin{lem}\label{Norm} With the previous notations, we have  canonical commutative diagram:
$$
\xymatrix{
  A_E\ar[d]^{\pi^E_{E'}}&\ar[l] \bigcap^d_{A_E} \hat\q_{E}^*\ar[r]^{\cap^q\mathbf c_E \qquad}\ar[d]\ar[d]& \colon \bigcap^d_{A_E}  H^1(E,ad(\rho_f))   \ar[d]^{\cap^q N^E_{E'}}\\
A_{E'}&\ar[l] \bigcap^d_{A_{E'}} \hat\q_{E'}^*\ar[r]^{\cap^q\mathbf c_{E'}\qquad}& H^1(E',ad(\rho_f))
}
$$
\end{lem}
\subsubsection{Compatible systems}
The goal of this section is to now explain the definition and construction of compatible systems under the norm map of zeta elements and Galois cohomology classes. We start by defining a class of set of abelian extensions of $F$.

\begin{defi}\label{S-admissible}
Let $S$ be a finite set of places of $F$.
A set $\mathfrak F$ of abelian extensions of $F$ is said $S$-admissible if
\begin{itemize}
\item[(i)]  For all $E\in\mathfrak F$, $E/F$ is unramified at places in $S$,
\item[(ii)] For any $E,E'\in\mathfrak F$ with $E\subset E'$, if the index of ramification of a place $w$ of $E$ is divisible by $p$, then $E'$ contains the maximal $p$-abelian extension of $F$ ramified at the place of $F$ dividing $w$.\end{itemize}
For $E'/E$ as in (ii), we denote by $Ram(E'/E)$ the set of places of $F$ which ramifies in $E'$ but not in $E$.

\end{defi}
In the next sections, we will discuss how one can choose elements $\xi$ satisfying the condition of Definition \ref{zeta-def} in an optimal way when $E$ varies in an $S$-admissible set of abelian extension of $F$. 
We will explain how this is connected
to the existence of an Euler system of rank $d$ whose definition we now recall.

\begin{defi}
Let $S$ be a finite set of primes of $F$ and  $(\rho_V,V)$ be a $p$-adic representation of $G_F$ which is unramified away from finitely many primes $S$ and $p$:
$$\rho_V\colon G_F\to GL_L(V)$$
 Let $\mathfrak F$ an $S$-admissible  set of abelian extensions of $F$.
An Euler system of rank $d$ for a  $G_F$-stable $O$-lattice  $T\subset V$ is a collection of elements $(z_E)_{E\in\mathfrak F}$
with $z_E\in\cap_{A_E}H^1(E,T)$ for each $E\in \mathfrak F$ such that for any $E',E\in\mathfrak F$ with $E'\supset E$, we have:
$$\cap^q_{A_E}  Cores_E^{E'} (z_{E'})= \prod_{v\in  Ram(E'/E)} P_v(\sigma_{v,E};V)\cdot \xi_E$$
where
\begin{itemize}
\item $\cap^q_{A_E}  Cores_E^{E'}$ is the norm map induced by the co-restriction map  in Galois cohomology
$$ \cap^q_{A_E} Cores^{E'}_E\colon H^1(E,T)\to \cap^q_{A_{E'}} H^1(E',T),$$
\item $Ram(E'/E)$ is the set of places of $E$ not dividing $p$ which ramifies in $E'/E$ and which are unramified in $E/F$,
\item $P_v(X;V)=det_V(1-X \rho_V(\Frob_v))$,
\item $\sigma_{v,E}$ is the geometric Frobenius in $Gal(E/F)$ at $v$.
\end{itemize}
\end{defi}

\begin{defi} \label{norm-comp-elements} Let $\mathfrak F$ as above.
A compatible system of  equivariant numbers for $V$, $S$ and $\mathfrak F$ 
is the data of  elements $\xi_E\in  A_E$ for each $E\in\mathfrak F$
such that for any $E',E\in\mathfrak F$ with $E'\supset E$, we have:
$$N_{E'/E}(\xi_{E'})= \prod_{v\in  Ram(E'/E)} P_v(\sigma_{v,E};V)\cdot \xi_E$$
where $N_{E'/E}$ is the natural projection map $A_{E'}\rightarrow A_E$,

For $V=ad(\rho_f)$, we say that such a system is a compatible system of congruence numbers (resp. of congruence anihilators) if $\xi_E\in Fitt_{A_E}(\wp_E/\wp_E^2)$ (resp. if $\xi_E$ 
annihilates $\wp_E/\wp_E^2$)  for each $E\in\mathfrak F$.
\end{defi}

An typical  example of compatible system of elements in $A_E$ as in the following definition is given by the equivariant $L$-values as defined in \S\ref{twisted-L}.

\begin{coro} \label{zeta}
If $(\xi_E)_{E\in\mathfrak F}$ is a compatible system of congruence numbers, then $(z_{f,E,\xi_E})_{E\in\mathfrak F}$ is an Euler system of rank $d$.
\end{coro}
\proof It is a straightforward consequence from the definitions, the construction and Lemma \ref{Norm}.\qed
For any $h_E=\sum_{w|p}h_w\in \bigoplus_{w|p} H^1(I_{E_w},O)^{G_{F_w}/I_{F_w}}$, let us denote by $\ell_{h_E}$ the $O$-linear form
on $ \Omega_{\Lambda_E/O}\otimes O\rightarrow O$  the isomorphisms of Lemma \ref{local-def}.(iii). The following corollary gives
many Euler systems with prescribed conditions at each place dividing $p$.

\begin{coro} \label{rankone}
If $(\xi_E)_{E\in\mathfrak F}$  is a compatible system of congruence annihilators, then 
for any $h:=(h_E)_E\in \displaystyle\varprojlim_{ E} \bigoplus_{w|p} H^1(I_{E_w},O)^{G_{F_w}/I_{F_w}}$, then the system of classes $(c^h_E)_{E\in\mathfrak F}$ defined by 
$$c^h_E:=  \mathbf c_E(\xi_E\cdot \ell_{h_E})\qquad\forall E\in\mathfrak F$$
defines an Euler system of rank $1$.
Moreover, for each place $w$ of $E$ dividing $p$ above the place $v$ of $F$, we have $res_w(c^h_E)\in H^1(I_{E_w},\CF^0_v)$ and its image  $res^0_w(c^h_E)$  
by the map $H^1(E_w,\CF_v^0)\rightarrow H^1(I_w,Gr^0_v)=H^1(I_w,O)$ is given by
$$res^0_w(c^h_E)=\xi_E\cdot h_w.$$
\end{coro}
\proof We just need to ensure the definition is meaningful. The norm relations on the classes will follow from the definitions and the diagram \eqref{Norm-diag}.
For each $E\in\mathfrak F$, we have an exact sequence
$$o\to \hat\q_E^*\to( \Omega_{\Lambda_E/O}\otimes O)^*\to (\wp_E/\wp_E^2)^\vee\to 0 $$
Since $\xi_E$ annihilates $\wp_E/\wp_E^2$, we deduce that
$\xi_E\cdot \ell_{h,E}\in \hat\q_E^*$ and therefore we can evaluate $\mathbf c_E$ on this elements.
The local property follows from the discussions of Corollary \ref{local-nontriv0} and the discussion preceeding it. \qed

\subsection{Variants}
\subsubsection{Iwasawa theory variant}
We now give a variant of the construction of the constructions made in the previous section. For any field totally real field $E$ as in the previous section, 
we denote by $E_\cyc$ the $\Zp$-cyclotomic extension of $E$, and by $E_n\subset E_\cyc$ the intermediate extension of degree $p^n$ over $E$.
We define the object $A_E^\Iw$, $\hat\q_E^\Iw$, $\q_E^\Iw$, $c_E^\Iw$, $H^1_\Iw(E, \cdot )$, the projective limit over $n$  for the norm or correstriction maps of the objects $A_{E_n}$,
$\hat\q_{E_n}$, $\q_{E_n}$, $c_{E_n}$, $H^1_\Iw(E_n, \cdot )$. It is straightforward to prove by taking projective limits to see that the results and constructions of the previous section go {\it mutatis mutandis}  with the obvious modifications.

\subsubsection{Hida family variant}
Similarly we can define the objects and construction of the previous sections and its Iwasawa theory variant  by replacing $f$ by  families of modular forms. We leave this to the interested reader.

\subsection{Deformations rings and Selmer groups}\label{kahler-diff}
The goal of this section is to recall the relation between the modules $\q_E$'s and Selmer groups in orders to exhibit
system of annihilators of congruence modules.

Let $LCN_O$ be the category of local noetherian complete $O$-algebra with  residue field at $\kappa=O/\varpi_O$.
Let $E/F$ as in the previous sections. We denote by $Ram(E/F)$ the places of $F$ which are ramified in the extension $E/F$.
Recall that $S$ is the finite set of places of $F$  dividing $\n p$. We consider the deformation functor:
$$\CF_E^S\colon LCN_O\longrightarrow Set$$
where $ \CF_E^S(A)$ is the set of strict equivalence classes of deformations
 $$\rho\colon G_E\rightarrow  GL_2(A)$$
such that
\begin{itemize}
\item $\rho\pmod{\m_A}\cong \rho_f \pmod {\varpi_O}$,
\item $\rho$ is unramified at places not dividing those in $S\cup\{p\}$,
\item $\rho$ is nearly ordinary at each place $w$ above $p$
\item $det\rho=\omega \epsilon_{cyc}$
\end{itemize}
It is well known by now classical arguments originally due Mazur, that this functor is representable by a local complete noetherian $O$-algebra $R_E$. See for example \cite{HidaSelmer3}. Let $\lb_E\colon R_E\to O$ be the homomorphism corresponding to the Galois representation $\rho_f$, and  $\wp_{R_E}:= Ker(\lb_E)$. It is proved in loc.cit.  using again arguments due to Mazur that there is an isomorphism
$$Hom_O(\wp_{R_E}/\wp_{R_E}^2,K/O)\cong Sel_S(E,ad(\rho_f))$$
where $Sel_S(E,ad(\rho_f)$ is the Selmer group defined as
 the kernel of the following restriction map:
\begin{eqnarray*}
H^1(Gal(\bQ/E), ad(\rho_f)\otimes K/O)&\rightarrow &\prod_{v\notin S_E} H^1(I_v, ad(\rho_f)\otimes K/O)\oplus \\ && \bigoplus_{v|p}H^1(I_v,(ad(\rho_v)/\CF_v^+(ad(\rho_v))\otimes K/O)
\end{eqnarray*}
For any extension $E/F$ which is unramified above places outside $S$ and those dividing $p$, we extend the definition  by defining $Sel_S(E,ad(\rho_f))$ to be the inductive limit of $Sel_S(E',ad(\rho_f))$ when $E$ runs in the set of finite extension of $F$ contained in $E$. In particular, we can define
$$Sel_S(E^\cyc,\rho_f):=\lim_{n} Sel_S(E_n,\rho_f)$$
It follows from result of Taylor-Wiles-Fujiwara and Hida \cite{HidaSelmer3}, that this Selmer group is co-torsion over $A_W^\cyc$. We denote by $\CL^{S,alg}_{E,f}$ the $A_E^\cyc$-Fitting ideal of the Pontryagin dual of $Sel_S(E^\cyc,\rho_f)$.

\begin{lem} \label{Fitt-Selmer}
Let $E$ and $S$ as above, then we have $\CL^{S,alg}_{E,f}\subset Fitt_{A_E^\cyc}(\q_E^\Iw)$.
\end{lem}
\proof
For each integer $n$, by the universal property of $R_{E_n}$, we have a surjective\footnote{The surjectivity is a classical fact obtained using the local  properties of the Galois representation
$\rho_{\T_{E_n}}$.} homomorphism:
$$ \psi_n \colon R_{E_n}\twoheadrightarrow \mathbf T_{E_n}$$
such that $\lb_{En}:=\lb_{f_{E_n}}\circ\psi_n$ .  By surjectivity of $\psi_n$ ,  we have a surjection of $A_{E_n}$-torsion modules:
\begin{eqnarray*}
\wp_{R_{E_n}}/\wp_{R_{E_n}}^2\twoheadrightarrow\wp_{f_{E_n}}/\wp_{f_{E_n}}^2=\q_{E_n}
\end{eqnarray*}
By passing to the projective limit over $n$, we deduce  a surjective homomorphism of $A_E^\cyc$-torsion modules:
$$Sel_S(E^\cyc,ad(\rho_f))^\vee\twoheadrightarrow \q_E^\Iw$$
which implies our claim.\qed

\subsection{Proof of the main Theorems}

\subsubsection{} Using the Iwasawa version of Corollary \ref{zeta} and \ref{rankone}, it just suffices to construct compatible system of congruence number in $A_E^\Iw$.
We will do so using  Lemma \ref{Fitt-Selmer}. We start by the following elementary lemma.
\begin{lem}\label{nonzerodiv}
For each $E\in \mathfrak F$ like in Definition \ref{norm-comp-elements}, let $I_E\subset A_E^\Iw$ be a non zero ideal such  that for any $E',E\in\mathfrak F$ such that $E\subset E'$, we have
$\pi_{E',E} (I_{E'})=Q_{E',E}I_{E}$  for some $Q_{E',E}\in A_{E}^\Iw$ satisfying
\begin{itemize}
\item $Q_{E',E}$ is not a zero divisor in $A_E^\Iw$
\item For any $E,E',E''\in\mathfrak F$ with $E\subset E'\subset E''$, we have 
$$Q_{E'',E}=\pi_{E'',E'}(Q_{E'',E'})Q_{E',E}.$$
\end{itemize}
Then there exists  a non trivial compatible system of elements $(\xi_E)_{E\in \mathfrak F}$ such that $\xi_E\in I_E$ for each $E\in\mathfrak F$ satisfying the relation
$\pi_{E',E} (\xi_{E'})=Q_{E',E}.\xi_{E}$ for any $E',E\in \mathfrak F$ such that $E\subset E'$, for any choice of $\xi_F\in I_F$.
\end{lem}
\proof For each $E',E$ as above let $\phi_{E',E}\colon I_E\to I_{E'}$ such that $\pi_{E',E}(x)=Q_{E',E}\phi_{E',E}(x)$ for any $x\in I_E$. This is well defined
since we assumed that the $Q_{E',E}$'s are not zero divisors. Moreover $\phi_{E',E}$ are surjective by hypothesis, and they define a projective system for $(I_E)_{E\in \mathfrak F}$ with surjective transition maps. Therefore we have a surjection$$\displaystyle\varprojlim_{E\in \mathfrak F, \phi_\cdot}I_E\twoheadrightarrow I_F$$ 
which implies our claim.\qed

\begin{prop}\label{selmerprop}
Let $\mathfrak F$ be a $S$-admissible set of abelian extension of $F$.
Let $E',E\in\mathfrak F$ with $E\subset E'$ as above.  Then we have 
$$\pi_{E',E}(\CL^{S,alg}_{E,f})= \prod_{v\in Ram(E'/E)}P_ v(q_v^{-1}\sigma_{v,E'}^\Iw) \cdot \CL^{S,alg}_{E',f}$$
where the product is on places $Ram(E'/E)$ which are ramified in $E'$ but not in $E$.
\begin{itemize}
\item $P_v(X,ad(\rho_f)):=det(1-X.ad(\rho_f)(Frob_v) )$,
\item $\sigma_{v,E'}^\Iw$ is the image of a geometric Frobenius automorphism at $v$ in $Gal(E'_\cyc/F)\subset (A_{E'}^\Iw)^\times$.
\end{itemize}
\end{prop}
\proof By induction on the cardinality of the set $Ram(E'/E)$, it is sufficient to prove the result in the following two cases.
\begin{itemize}
\item[(a)] $Ram(E'/E)=\emptyset$,
\item[(b)]  $ Ram(E'/E)=\{v_0\}$,
\end{itemize}
The case (a) follows from the isomorphism below which is valid for any finite set $S$ containing the places dividing the level of $f$
$$ Sel_S(E'_\cyc,ad(\rho_f))= Sel_S(E_\cyc,ad(\rho_f))^{Gal(E'/E)} $$
which follows easily from the inflation-restriction exact sequence since    $ad(\bar\rho_f)$ has no invariants by $G_{E'}$, and because 
the ramification conditions are the same on both Selmer groups since $E'/F$ is unramified at all places where $E/F$ is (see \cite{HidaSelmer3} for example). 
Moreover a class in $H^1_\Iw(E,\rho_f)$ which is unramified at a place $w'$ after restriction to $E'$ is already unramified at the place $w$ of $E$ below $w'$ since the ramification index is prime to $p$, by the hypothesis of the set $\mathfrak F$.

For the case (b), we may assume that $E'= EF_0$ where $F_0$ is the $p$-maximal abelian extension of $F$ unramified away from $v_0$. In that case, we have
$$  Sel_{S\cup\{v_0\} }(E_\cyc,ad(\rho_f))= Sel_S(E'_\cyc,ad(\rho_f))^{Gal(E'/E)} $$
We argue as in the isomorphism used for (a). We only need to show that the local condition of a class in $H^1_\Iw(E,ad(\rho_f)$ that is ramified possibly at places dividing $v_0$ becomes unramified after restrictipn to $E'$. This is obvious since $E'$ contains the maximal $p$-abelian extension of $F$ that is unramified away from $v_0$. From the previous ismorphism, we deduce that
\begin{equation}\label{equa1}
\pi_{E',E}(\CL^{S,alg}_{E,f})= \CL^{S\cup\{v_0\},alg}_{E',f}
\end{equation}
Now consider the exact sequence
$$
0 \to Sel_{S }(E_\cyc,ad(\rho_f)) \to Sel_{S\cup\{v_0\} }(E_\cyc,ad(\rho_f)) \to H^1(I_{v_0}, ad(\rho_f)\otimes (A_E^\Iw)^\vee)^{G_{v_0}/I_{v_0}}\to 0
$$
The left exactness follows from the definition of the ramification conditions of the Selmer group on the left and the fact due to Shapiro's lemma that
$$\varinjlim_{ n} H^1(E_n,ad(\rho_f)\otimes K/O)= H^1(F, ad(\rho_f)\otimes_O(A_E^\Iw)^\vee)$$
and
$$\bigoplus_{w|v_0}(H^1(I_w,ad(\rho_f)\otimes \Lambda^\vee)^{D_v/I_v})^\vee
=H^1(I_{v_0} ad(\rho_f)\otimes O[Gal(F_\m/F)]\otimes \Lambda^\vee)^{D_\q/I_\q})^\vee
$$
where $w$ runs in the sate of places of $E$ above $v_0$. The exactness on the right is a Theorem of Greenberg-Vatsal  \cite{GV} see also section 4.2 of \ref{Greenberg}, from which we get 
\begin{equation}\label{equa2}
\CL^{S\cup\{ v_0 \},alg}_{E,f}=P_ v(q_{v_0}^{-1}\sigma_{v_0,E'}^\Iw) \cdot \CL^{S,alg}_{E,f}.
\end{equation}
Combining \eqref{equa1} and \eqref{equa2}, we get (b).\qed

Since the Fitting ideal of the right hand side of the equality above 
is given by $P_{\mathfrak q}(N\mathfrak q^{-1}<\mathfrak q>_{\m})$ which follows from  a classical calculation (see for example  \cite{GV}),
the claim  is now  a direct consequence of the assertions (ii) and (iii) of Proposition \ref{selmerprop}.
\qed

\subsubsection{Proof of Theorem \ref{thm2} and \ref{thm3}}
Note that Lemma \ref{nonzerodiv} applies with $I_E= \CL^{S,alg}_{E,f}.$ as  it is straightforward to see that $P_v(q_v^{-1}\sigma_{v,E})$ is not a zero divisor in $A_E^\Iw$ since its image on every irreducible component of $A_E^\Iw$ is clearly non-zero.
By Proposition \ref{selmerprop} and Lemma \ref{nonzerodiv}, there exists therefore 
a compatible system of congruence numbers $(\xi_E)_{E\in\mathfrak F}$ such that $\xi_F=\CL_{f}^{S,alg}\subset \Lambda= A_F^\Iw=O[[Gal(F_\cyc/F)]]$.
It now suffices to apply  the Iwasawa variant of Corollaries \ref{zeta} and \ref{rankone} to this system, and we are done.\qed

\subsubsection{Proof of Theorem \ref{thm4} and \ref{thm5}}

The conjectures \ref{equivariant-conj} and \ref{equivariant-conj-weak} state that $\CL_{E,f}^{S,an}$ is  respectively a compatible system
of congruence numbers and of congruence annihilators.
It therefore  suffices to apply the Iwasawa variant of Corollaries \ref{zeta} and \ref{rankone} to this system.\qed

\end{document}